\documentclass[11pt,aip]{revtex4-1}

\usepackage{bm, amsfonts, mathtools}
\usepackage{graphicx, color, wasysym}
\usepackage{textcomp}
\usepackage{amsmath, amssymb, amsthm, latexsym, epsfig}
\usepackage{appendix}
\usepackage[normalem]{ulem}
\usepackage{soul,xcolor}

\DeclareMathOperator{\D}{d\!}
\DeclareMathOperator{\E}{e} 
\DeclareMathOperator{\I}{i}

\begin{document} 

\title[Can Umbral and $q$-calculus be merged?]{Can Umbral and $q$-calculus be merged?}

\author{G. Dattoli} 

\affiliation
{ENEA - Centro Ricerche Frascati, via E. Fermi, 45, IT 00044 Frascati (Roma), Italy \vspace{2mm}}

\author{B. Germano}
\author{M. R. Martinelli}
\affiliation
{University of Rome, La Sapienza, Department of Methods and Mathematic Models for Applied Sciences, Via A. Scarpa, 14, 00161 Rome, Italy \vspace{2mm}}

\author{K. G\'{o}rska}
\affiliation
{H. Niewodnicza\'{n}ski Institute of Nuclear Physics, Polish Academy of Sciences, 
ul.Eljasza-Radzikowskiego 152, PL 31342 Krak\'{o}w, Poland \vspace{2mm}}

\email{giuseppe.dattoli@enea.it; bruna.germano@sbai.uniroma1.it; katarzyna.gorska@ifj.edu.pl; martinelli@dmmm.uniroma1.it}

\begin{abstract}
The $q$-calculus is reformulated in terms of the umbral calculus and of the associated operational formalism. We show that new and interesting elements emerge from such a restyling. The proposed technique is applied to a different formulations of $q$ special functions, to the derivation of integrals involving ordinary and $q$-functions and to the study of $q$-special functions and polynomials.
\end{abstract}
%
\maketitle

\section{Introduction}

In this paper we explore the possibility of extending the umbral point of view to the $q$-calculus \cite{FHJacson04}. The results we obtain are fairly encouraging and we consider worth to pursue the relevant exposition. 
The $q$-numbers and their associated factorial are defined below
\begin{equation}\label{eq1}
[n]_{q} = \frac{1-q^{n}}{1-q} \qquad \text{and} \qquad [n]_{q}! = \prod_{r=1}^{n} [r]_{q}, \quad \text{where} \quad 0 < q < 1.
\end{equation}
A previous investigation in this direction was addressed by Roman in Ref. \cite{SRoman-b-84, SRoman85}, here we develop an alternative point of view, based on a different conception of umbral calculus \cite{SLicciardi18}.

In the following, we exploit a function, namely the $q$-Bessel like function, which provides the pivot tool to establish the relevance of umbral methods to $q$-calculus. We start with studying the simplest case of $q$-Bessel type function written as
\begin{equation}\label{eq3}
C_{0}^{(q, 1)}(x) = \sum_{r=0}^{\infty} \frac{(-x)^r}{r! [r]_{q}!}, \qquad |x| < \infty,
\end{equation}
which is referred also as $(q,1)$-Tricomi function. 
In the language of the umbral calculus we define the umbral image of the function in Eq. \eqref{eq3} as follows \cite{GDattoli15, DBabusci13}
\begin{equation}\label{25/03-1}
C_{0}^{(q, 1)}(x) = \E^{-{_{q}\hat{c}_{z}} x} \varphi_{0} = \sum_{r=0}^{\infty} \frac{(-x)^{r}}{r!}\, {_{q}\hat{c}_{z}}^{\!r} \varphi_{0}
\end{equation} 
with ${_{q}\hat{c}_{z}}$ being an operator defined below. 
The series expansion on the RHS (right hand site) of the above function, is the crucial step of our discussion. In Eq. \eqref{25/03-1}, as 
developed in previous researches (see \cite{DBabusci17} and references therein), we can say that the $q$-operator ${_{q}\hat{c}}_{z}$ acts on the vacuum  $\varphi_{0}$ according to the rule ${_{q}\hat{c}_{z}}^{\!\mu} f(z)\vert_{z=0} = {_{q}\hat{c}_{z}}^{\!\mu} \varphi_{0} = f(\mu)$ which for integer $m$ and $f(z) = ([z]_{q}!)^{-1}$ gives
\begin{equation}\label{eq5a}
{_{q}\hat{c}_{z}}^{\!m} \varphi_{0} = \frac{1}{[m]_{q}!}. 
\end{equation}
In general the exponent $\mu$ is not restricted to integers but it can be any real. 
The extension of the identity \eqref{eq5a} occurs through the use of the Thomae-Jackson $q$-Gamma function \cite{FHJacson04, JThomae1869}, namely
\begin{equation}\label{eq5b}
{_{q}\hat{c}_{z}}^{\mu} \varphi_{0} = \frac{1}{{_q\Gamma}(1 + \mu)} , \qquad {_q\Gamma}(x) = (1-q)^{1-x}\prod_{n=0}^{\infty} \frac{1 - q^{n+1}}{1 - q^{n+x}}.
\end{equation}
The usefulness of the umbral image restyling stays in the fact that we have associated to the $q$-function a corresponding "image" with known properties. 
We can therefore take advantage from such a correspondence to infer the properties of $q$-function under study by translating to the $q$-functions those characterizing the associated images. Taking notice of the identity in Eq. \eqref{25/03-1} we recover the series \eqref{eq3} from the corresponding umbral images, whose Taylor expansion around $x=0$ yields 
\begin{equation}\label{25/03-2}
\E^{- {_{q}\hat{c}_{z}} x} \varphi_{0} = \sum_{r=0}^{\infty} \frac{(-x)^{r}}{r!}  {_{q}\hat{c}_{z}}^{r} \varphi_{0} = \sum_{r=0}^{\infty} \frac{(-x)^{r}}{r! [r]_{q}!}.
\end{equation}
Moreover, calculating the Borel transform of the $(q, 1)$-Tricomi function:
\begin{align}\label{25/03-3}
\begin{split}
\int_{0}^{\infty} \E^{-s} C_{0}^{(q, 1)}(x s) \D s = \int_{0}^{\infty} \E^{-s} \E^{- {_{q}\hat{c}_{z}} x s} \varphi_{0} \D s = \int_{0}^{\infty} \E^{-s - {_{q}\hat{c}_{z}} x s} \D s \varphi_{0} \\ = \frac{1}{1 + {_{q}\hat{c}_{z}} x} \varphi_{0} = {_{q}e}(x)
\end{split}
\end{align}
we introduce the new function $ {_{q}e}(x)$ whose umbral representation is given by the second line of Eq. \eqref{25/03-3}.
The Taylor expansion of ${_{q}e}(x)$ around $x=0$ reads
 \begin{equation}\label{25/03-4}
 \frac{1}{1 + {_{q}\hat{c}_{z}} x} \varphi_{0} = \sum_{r=0}^{\infty} (-x)^{r} {_{q}\hat{c}_{z}}^{r} \varphi_{0} = \sum_{r=0}^{\infty} \frac{(-x)^{r}}{[r]_{q}!},
 \end{equation}
that gives the series representation of the $q$-exponential converging provided that $|x| < (1+q)^{-1}$. The Borel transform given by Eq. \eqref{25/03-3} can be viewed as the integral representation of ${_{q}e}(x)$ function. The advantage of the integral representation, with respect to the series expansion, stems from the fact that no convergence problems arise. In the following we use, for numerical computations involving the  $q$-exponential, the relevant integral representation. The series expansion will be exploited, in the quoted range of convergence, for analytical computations.  \\
\noindent
The chain of identities leading to Eq. \eqref{25/03-3} can be summarized as it follows:
\begin{itemize}
\item[(1)] we have used the umbral form of the $(q, 1)$-Tricomi function; 
\item[(2)] we have assumed that the vacuum can be brought outside the integral;
\item[(3)] we have used the Laplace transform; 
\item[(4)] we obtained the umbral representation of ${_{q}e}(x)$.
\end{itemize}
For latter convenience we premise the following identity from the elementary calculus.

We have associated to the exponential $q$-function a rational image function, which can be helpful for the study of the relevant properties. The Lorentzian function can for example be exploited as image of the $q$-Gaussian ${_{q}e}(x^{2}) = (1 + {_{q}\hat{c}_{z}}x^{2})^{-1}\varphi_{0}$. A paradigmatic use of the image correspondence is realized through the following example noting the identity 
\begin{equation}\label{eq8}
I_{L}(a) = \int_{-\infty}^{\infty} \frac{\D x}{1 + a x^{2}}  = \frac{\pi}{\sqrt{a}}, \qquad a > 0.
\end{equation}
If we are interested in evaluating the infinite integral of ${_{q}e(x^{2})}$ we find
\begin{equation}\label{eq9}
\int_{-\infty}^{\infty} {_{q}e(x^{2})} \D x = \int_{-\infty}^{\infty} \frac{\D x}{1 + {_{q}\hat{c}_{z}} x^{2}}\varphi_{0} = I_{L}({_{q}\hat{c}}) \varphi_{0} = \pi\, {_{q}\hat{c}_{z}}^{-1/2} \varphi_{0} = \frac{\pi}{{_{q}\Gamma}(1/2)}.
\end{equation}
We underscore once more that the key note of the operation has been treating the operator ${_{q}\hat{c}_{z}}$ as a constant. In the case of the ordinary Gamma function  $\Gamma(1/2) = \sqrt{\pi}$, we assume that something analogous does hold for the $q$-Gamma extension and set
\begin{equation}\label{eq10}
{_{q}\Gamma(1/2)} = \sqrt{\pi_{q}}, \qquad \lim_{q\to 1^{-}} {_{q}\Gamma(1/2)} = \sqrt{\pi}.
\end{equation}
The "protocol" for the evaluation of $\sqrt{\pi_{q}}$ occurs as noted below. According to the second of Eqs. \eqref{eq5b} we write
\begin{align}\label{eq11}
\begin{split}
{_{q}\Gamma}(1/2) & = (1-q)^{1-1/2} \prod_{n=0}^{\infty} \frac{1-q^{n+1}}{1-q^{n+1/2}} = \sqrt{1-q} \prod_{n=0}^{\infty} \frac{1-\sqrt{q}^{(2n + 2)}}{1 - \sqrt{q}^{(2n+1)}} \\ &= \sqrt{1-q} \prod_{n=0}^{\infty} \frac{[2n + 2]_{\sqrt{q}}}{[2n + 1]_{\sqrt{q}}}.
\end{split}
\end{align}
The conclusion of the previous discussion is that
\begin{equation}\label{eq12}
\int_{-\infty}^{\infty} {_{q}e(x^{2})} \D x = \frac{\pi}{\sqrt{\pi_{q}}},
\end{equation}
whose limit correctly yields (see Fig. \ref{fig1} where we have reported $\pi_{q}$ in the interval $q = 0, 1$)
\begin{equation}
 \lim_{q\to 1^{-}} \int_{-\infty}^{\infty} {_{q}e(x^{2})} \D x = \lim_{q\to 1^{-}} \frac{\pi}{\sqrt{\pi_{q}}} = \sqrt{\pi}
\end{equation}
The "experimental" check of the correctness of the result in Eq. \eqref{eq12} is obtained by computing ${_{q}\Gamma(1/2)}$ by means of Eq. \eqref{eq5b} and then by confronting the numerical value from Eq. \eqref{eq11} with that deriving from a direct numerical integration.

The next step is that of checking whether the procedure we have envisaged is consistent with the derivation of other integrals involving $q$-functions whose ordinary counterpart contains $\sqrt{\pi}$. We will indeed establish the "complete" $q$-Fresnel integrals
\begin{equation}\label{eq13}
{_{q}I^{cos}_{F}} = \lim_{x\to\infty} \int_{0}^{x} {_{q}\cos}(y^{2}) \D y = \frac{1}{2} \int_{-\infty}^{\infty}  {_{q}\cos}(y^{2}) \D y
\end{equation}
and
\begin{equation}\label{eq13a}
{_{q}I^{sin}_{F}} = \lim_{x\to\infty} \int_{0}^{x} {_{q}\sin}(y^{2}) \D y = \frac{1}{2} \int_{-\infty}^{\infty}  {_{q}\sin}(y^{2}) \D y.
\end{equation}
The proof is achieved by noting that
\begin{equation}\label{eq14}
{_{q}\cos(x)} = \frac{1}{1 + x^{2}\, {_{q}\hat{c}_{z}}^{2}} \varphi_{0} \qquad \text{and} \qquad {_{q}\sin(x)} = \frac{x\, {_{q}\hat{c}_{z}}}{1 + x^{2}\, {_{q}\hat{c}_{z}}^{2}} \varphi_{0} 
\end{equation}
which satisfy like to standard exponential function the property ${_{q}\cos}(x) + \I {_{q}\sin}(x) = {_{q}e}(-\I\!x)$, where the umbral representation of ${_{q}e}(x)$ is given by Eq. \eqref{25/03-3}. The associated integrals writes
\begin{equation}\nonumber
\int_{-\infty}^{\infty}\!\! {_{q}\cos}(y^{2}) \D y = \int_{-\infty}^{\infty} \frac{\D y}{1 +  {_{q}\hat{c}_{z}}^{2} y^{4}} \varphi_{0} = \int_{0}^{\infty}\!\!\!  \E^{-s} \left[\int_{-\infty}^{\infty}\!\!\! \E^{-s\, {_{q}\hat{c}_{z}}^{2} y^{4}} \D y\right]\! \D s\, \varphi_{0}  = \frac{\Gamma(1/4) \Gamma(3/4)}{2\,{_{q}\Gamma}(1/2)} 
\end{equation}
and
\begin{equation}\nonumber
\int_{-\infty}^{\infty}\!\! {_{q}\sin}(y^{2}) \D y = \int_{-\infty}^{\infty} \frac{y^{2} {_{q}\hat{c}_{z}} \, \D y}{1 +  {_{q}\hat{c}_{z}}^{2} y^{4}} \varphi_{0} = \int_{0}^{\infty}\!\!\!  \E^{-s} \left[{_{q}\hat{c}_{z}} \int_{-\infty}^{\infty}\!\!\! y^{2} \E^{-s\, {_{q}\hat{c}_{z}}^{2} y^{4}} \D y\right]\! \D s\, \varphi_{0} = \frac{\Gamma(1/4) \Gamma(3/4)}{2\,{_{q}\Gamma}(1/2)}.
\end{equation}
The use of the identity 
\begin{equation}\label{eq16}
\Gamma(x) \Gamma(1-x) = \pi/\sin(\pi x)
\end{equation} 
for $x = 1/4$ yields to ${_{q}I^{cos}_{L}} = {_{q}I^{sin}_{L}} = \pi/\sqrt{8 \pi_{q}}$ and {for $q\to 1^{-}$ we have a direct extension of the Fresnel integrals} ${_{1^{-}}I^{cos}_{L}} = {_{1^{-}}I^{sin}_{L}}= \sqrt{\pi/8}$. 
A further result obtained by a straightforward application of the method is
\begin{equation}\label{eq17}
\int_{0}^{\infty} {_{q}e}(x^{m}) \D x = \frac{1}{m}\, {_{q}\Gamma}\left(\frac{m-1}{m}\right).
\end{equation}

The next step is the derivation of the integral 
\begin{equation}\label{eq19}
I_{0}^{(q, 1)} = \int_{-\infty}^{\infty} C_{0}^{(q, 1)}(x^{2}) \D x,
\end{equation}
which requires as premise the Gaussian integral
\begin{equation}\label{19a}
I_{G}(a) = \int_{-\infty}^{\infty} \E^{-a x^{2}} \D x = \sqrt{\frac{\pi}{a}}, \qquad a > 0,
\end{equation}
accordingly we find
\begin{equation}\label{eq19b}
I_{0}^{(q, 1)} = \int_{-\infty}^{\infty}  {_{q}e}^{- {_{q}\hat{c}_{z}} x^{2}} \D x\, \varphi_{0} = \sqrt{\frac{\pi}{_{q} \hat{c}_{z}}}\, \varphi_{0} = \sqrt{\frac{\pi}{\pi_{q}}}.
\end{equation}

Note that the use of the previous result, along with the integral in Eq. \eqref{25/03-3}, can be used to prove Eq. \eqref{eq12}, independently of the use of the Lorentz function as umbral image. We find indeed
\begin{align}\nonumber
\begin{split}
\int_{-\infty}^{\infty} {_{q}e}(x^{2}) \D x & = \int_{0}^{\infty} \E^{-s} \left[\int_{-\infty}^{\infty} C_{0}^{(q, 1)}(x^{2} s) \D x\right] \D s \\ 
& = \sqrt{\frac{\pi}{\pi_{q}}} \int_{0}^{\infty} s^{-1/2} \E^{-s} \D s = \frac{\pi}{\sqrt{\pi_{q}}}.
\end{split}
\end{align}
Most of the results obtained so far have $\pi_{q}$ as pivot \cite{WSChung14, RWGosper99}. We can therefore summarize  the discussion of this introductory section with a question about the meaning of $\pi_{q}$ and whether any geometrical interpretation (in trigonometric sense) can be appended to it. We will try to speculate on these points in the forthcoming parts of the paper.

Before concluding this section we consider the use of the umbral method for the evaluation of the successive derivatives of the  $q$-Gaussians, which opens further interesting possibilities, which will be exploited in the final section within a specific application. We take into account that \cite{GDattoli17a, dzielo}
\begin{equation}\label{3/04-1}
(-1)^{n} \partial_{x}^{n} \E^{-a x^{2}} = \E^{-a x^{2}} H_{n}(2ax, -a),
\end{equation}
where $H_{n}(x, y)$ is the Hermite polynomials in two variables \cite{dzielo} 
\begin{equation}\label{3/04-2}
H_{n}(x, y) = n! \sum_{r=0}^{\lfloor n/2 \rfloor} \frac{x^{n-2r} y^{r}}{(n-2r)! r!}.
\end{equation}
We apply the same formalism to get the successive derivatives of the $q$-Bessel type function, thus getting
\begin{align}\label{19/04-1}
\begin{split}
(-1)^{n}\, \partial_{x}^{n} C^{(q, 1)}(x^{2}) &= (-1)^{n} \partial_{x}^{\,n} \E^{- {_{q}c_{z}} x^{2}}\, \varphi_{0} \\
& = \E^{- {_{q}c_{z}} x^{2}} H_{n}(2 \, {_{q}c_{z}} x, - {_{q}c_{z}})\, \varphi_{0}.
\end{split}
\end{align}
Using Eq. \eqref{3/04-2} for the Hermite polynomials we can rewrite Eq. \eqref{19/04-1} as
\begin{align}\label{23/04-1}
\begin{split}
(-1)^{n}\, \partial_{x}^{n}C^{(q, 1)}(x) & = n! \sum_{r=0}^{\lfloor n/2 \rfloor} \frac{(-1)^{r}\, (2x)^{n-2r}}{(n - 2r)! r!} {_{q}c_{z}}^{n-r} \E^{-{_{q}c_{z}} x^{2}} \varphi_{0} \\
& =  x^{-n/2} n! \sum_{r=0}^{\lfloor n/2 \rfloor} \frac{(-1)^{r}\, (2\sqrt{x})^{n-2r}}{(n - 2r)! r!} J_{n-r}^{(q, 1)}(2 x),
\end{split}
\end{align}
where the $q$-Bessel function has the below umbral form 
\begin{align}\label{3/04-6}
\begin{split}
{_{q}\hat{c}_{z}}^{\mu} \E^{-a {_{q}\hat{c}_{z}}} \varphi_{0} & = \sum_{r=0}^{\infty} \frac{(-a)^{r}}{r!} {_{q}\hat{c}_{z}}^{r + \mu} \varphi_{0} = \sum_{r=0}^{\infty} \frac{(-a)^{r}}{r! [\mu+r]_{q}!} \\
& = a^{-\mu/2} J_{\mu}^{(q, 1)}(2\sqrt{a}),
\end{split}
\end{align}
which for $\mu=0$ reduces to the $(q, 1)$-Tricomi function given via Eq. \eqref{25/03-2}.

\section{$q$-Wallis formula and $q$-sine infinite product formula}

In this section we go back to the definition of $\pi_{q}$, which represents an interesting key note to extend the role of $q-$calculus within its umbral perspective. The Wallis product is a well-known result of pre-calculus mathematics yielding $\pi$ in terms of an infinite product, namely 
\begin{equation}\label{20/03-1}
\frac{\pi}{2} = \prod_{n=0}^{\infty} \frac{(2n+2)^{2}}{(2n + 1) (2n + 3)} = \frac{2\cdot 2}{1\cdot 3} \cdot \frac{4 \cdot 4}{3 \cdot 5} \cdot \frac{6\cdot 6}{5 \cdot 7} \cdots.
\end{equation}
An analogous expression holds for the $q$-counterpart. The relevant proof is straightforward, noting indeed that
\begin{equation}\label{eq20}
1-q = (1 - \sqrt{q})\, \frac{1 - (\sqrt{q})^{2}}{1 - \sqrt{q}} = \frac{[2]_{\sqrt{q}}}{[\infty]_{\sqrt{q}}}, \qquad \text{where} \quad [\infty]_{q^{k}} = (1-q^{k})^{-1},
\end{equation}
and squaring Eq. \eqref{eq11} we find for $\pi_{q}$ the following infinite product identity
\begin{equation}\label{eq21}
\pi_{q} = \frac{[2]_{\sqrt{q}}}{[\infty]_{\sqrt{q}}} \prod_{n=0}^{\infty} \left(\frac{[2n + 2]_{\sqrt{q}}}{[2n + 1]_{\sqrt{q}}}\right)^{2},
\end{equation}
which is equivalent (albeit written in a different way) to an analogous expression derived by Chung, Kim and Mansour in Ref. \cite{WSChung14}.

Let us now try to get an answer whether we can push the analogies by getting a more profound meaning for $\pi_{q}$. Going back to the $q$-Gamma definition (Eq. \eqref{eq5b}) we note that it can be cast in the form
\begin{equation}\label{eq22}
\Gamma_{q}(x) = (1-q)^{1-x} \prod_{n=0}^{\infty} \left(\left[1 + \frac{x-1}{n+1}\right]_{q^{n+1}}\right)^{\!-1}.
\end{equation}
Using Eq. \eqref{eq16} as reference identity we define the sine-$q$ like function
\begin{equation}\label{eq23}
\sin_{q}(\pi_{q} x) = \frac{\pi_{q}}{\Gamma_{q}(x) \Gamma_{q}(1-x)} = \frac{\pi_{q}}{1-q} \prod_{n=1}^{\infty} \left(\left[1 + \frac{x-1}{n}\right]_{q^{n}}\, \left[1 - \frac{x}{n}\right]_{q^{n}}\right),
\end{equation}
which 
{can be expressed by the Gosper sin-$q$ function denoted here as ${^{G}\!\sin}_{q}(\pi x)$ \cite{IMezo13, RWGosper99} in the following way $\sin_{q^{2}}(\pi_{q} x) = {^{G}\!\sin_{q}(\pi x)}/q^{(x-1/2)^{2}}$}. It is indeed written in a form analogous to the Euler infinite product \cite{WFEberlein77} for the ordinary sine function. The relevant zeros are located at
\begin{equation}\label{eq24}
\sin_{q}(k \pi_{q}) = 0 \quad \text{and} \quad \sin_{q}[(2k + 1) \pi_{q}/2] = 1,  \qquad k = 0, 1, 2, \ldots.
\end{equation}
The function is by no means periodic, it has increasing oscillating amplitude (see Fig. \ref{fig1}) and the maxima are located at 
\begin{equation}\label{eq25}
y^{\star} = \frac{2 k + 1}{2} \pi_{q}, \qquad k = 0, 1, 2, \ldots.
\end{equation}
Using the same criterion (namely by following the analogy with the ordinary case) we can define the cos-$q$ analog, which writes
\begin{equation}\label{eq26}
\cos_{q}(\pi_{q} x) = \prod_{n=1}^{\infty} \left(\left[1 + \frac{x}{n - 1/2}\right]_{q^{n-1/2}} \left[1 - \frac{x}{n -1/2}\right]_{q^{n-1/2}}\right),
\end{equation}
{and gives $\cos_{q}(\pi_{q} x) = \pi_{q}/[\Gamma_{q}(1/2 + x) \Gamma(1/2-x)]$. The behavior {of Eq. \eqref{eq26}} is reported in Fig. \ref{fig1}(b), and discloses characteristics similar to those of $\sin_{q}(\pi_{q} x)$ given by Eq. \eqref{eq23} also reported in Fig. \ref{fig1}(a). The relation between $\cos_{q}(\pi_{q} x)$ and the Gosper cos-$q$ function ${^{G}\!\cos_{q}(\pi x)}$ \cite{RWGosper99}, is equal to $\cos_{q^{2}}(\pi_{q^{2}} x) = {^{G}\!\cos_{q}(\pi x)}/q^{x^{2}}$. The following properties are  worth to be mentioned
\begin{equation}\label{eq27}
\cos_{q}(k \pi_{q}/2) = 0, \quad k=1, 3, \ldots, 2m+1; \qquad \cos_{q}(0) = 1.
\end{equation}
It is furthermore important to stress that both functions have not definite parity properties under variable reflection $x\to -x$.
\begin{figure}[!h]
\begin{center}
\includegraphics[scale=0.50]{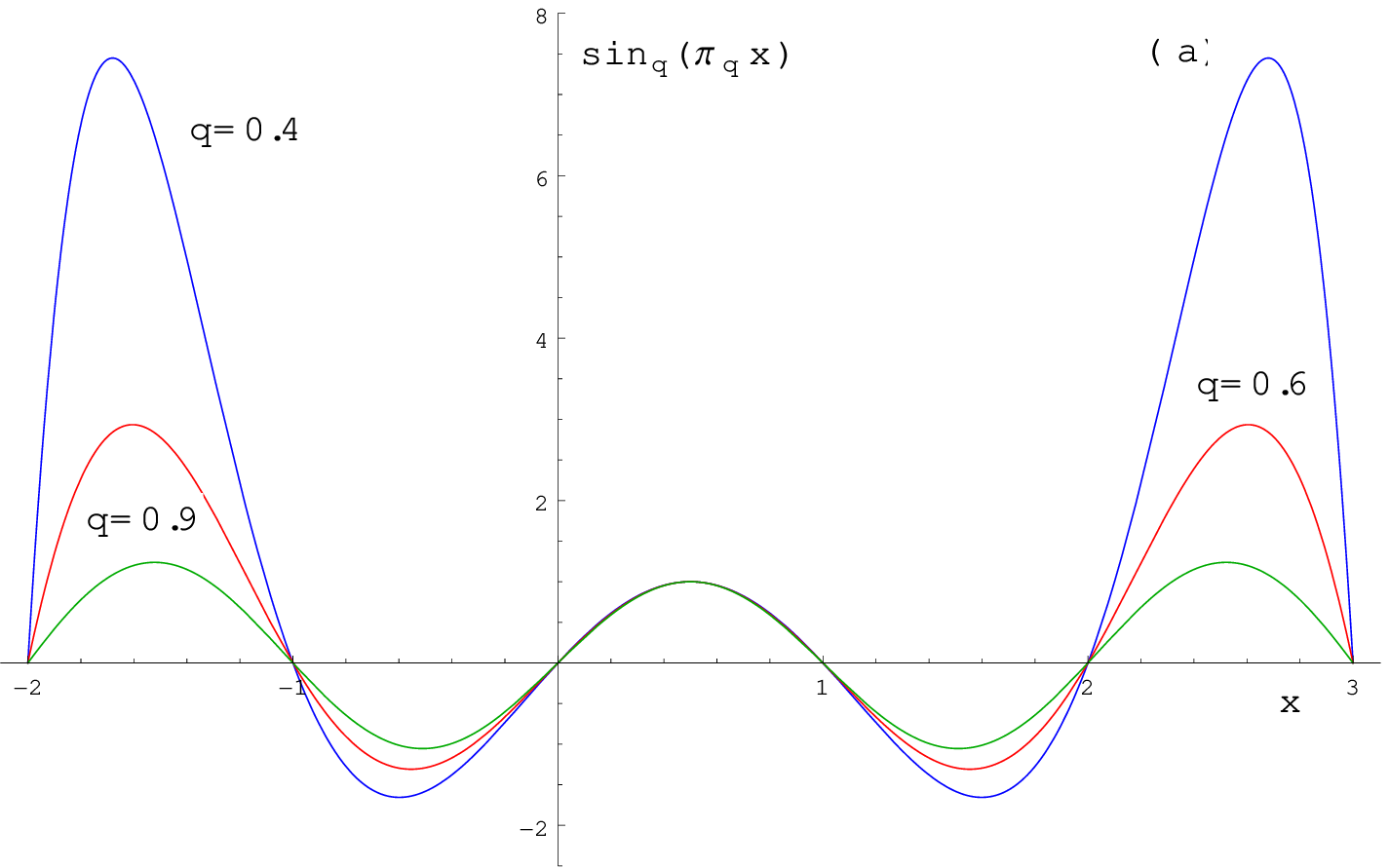}
\includegraphics[scale=0.48]{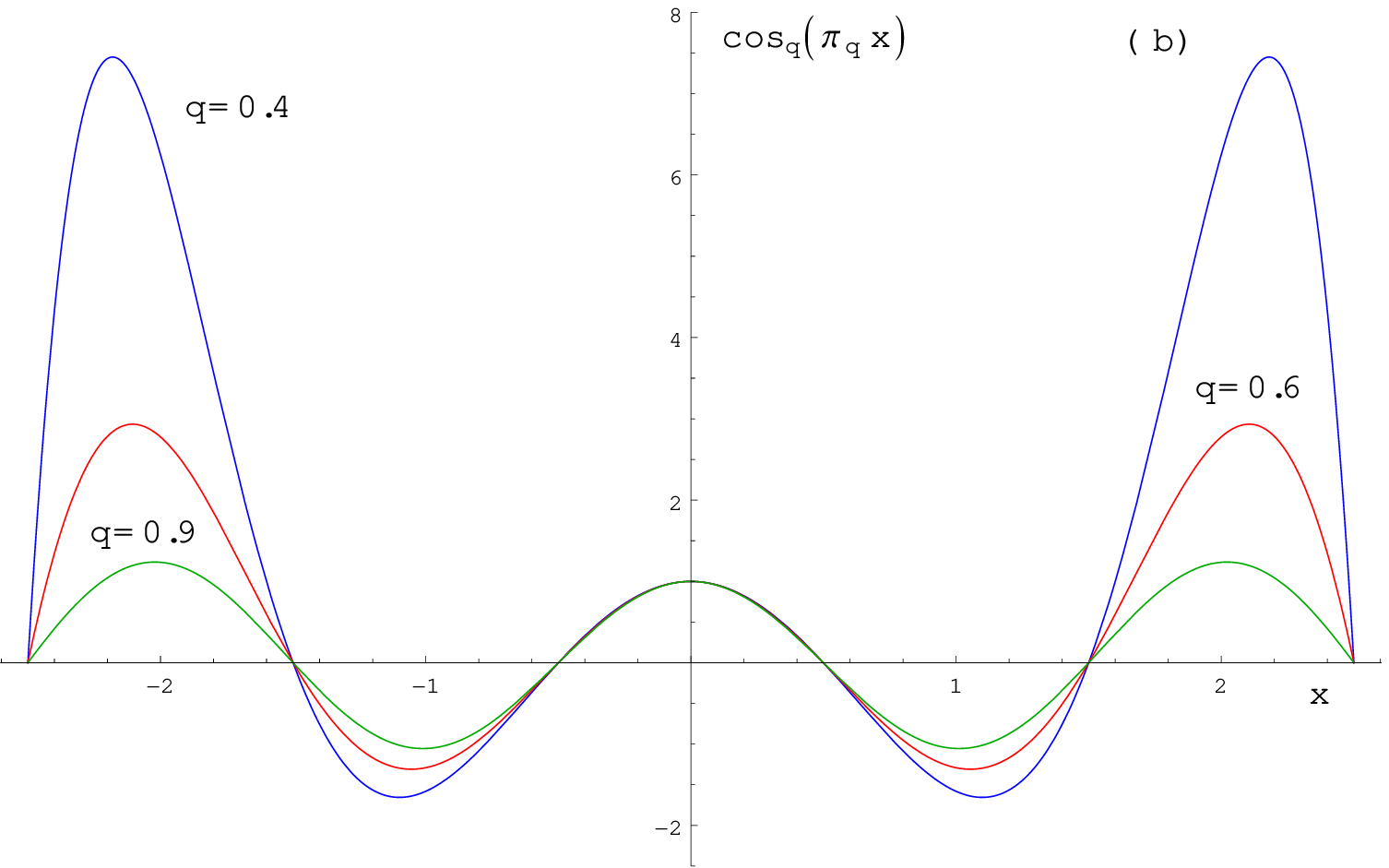}
\end{center}
\caption{\label{fig1} The plot of sine-$q$ (Fig. \ref{fig1}(a)) and cosine-$q$ (Fig. \ref{fig1}(b)) functions for $q=0.4$ (the blue curve), $q=0.6$ (the red curve), and $q=0.9$ (the green curve). It can be seen that $\sin_{q}(\pi_{q} x + 1/2) = \cos(\pi_{q} x)$.}
\end{figure}
\begin{figure}[!h]
\begin{center}
\includegraphics[scale=0.60]{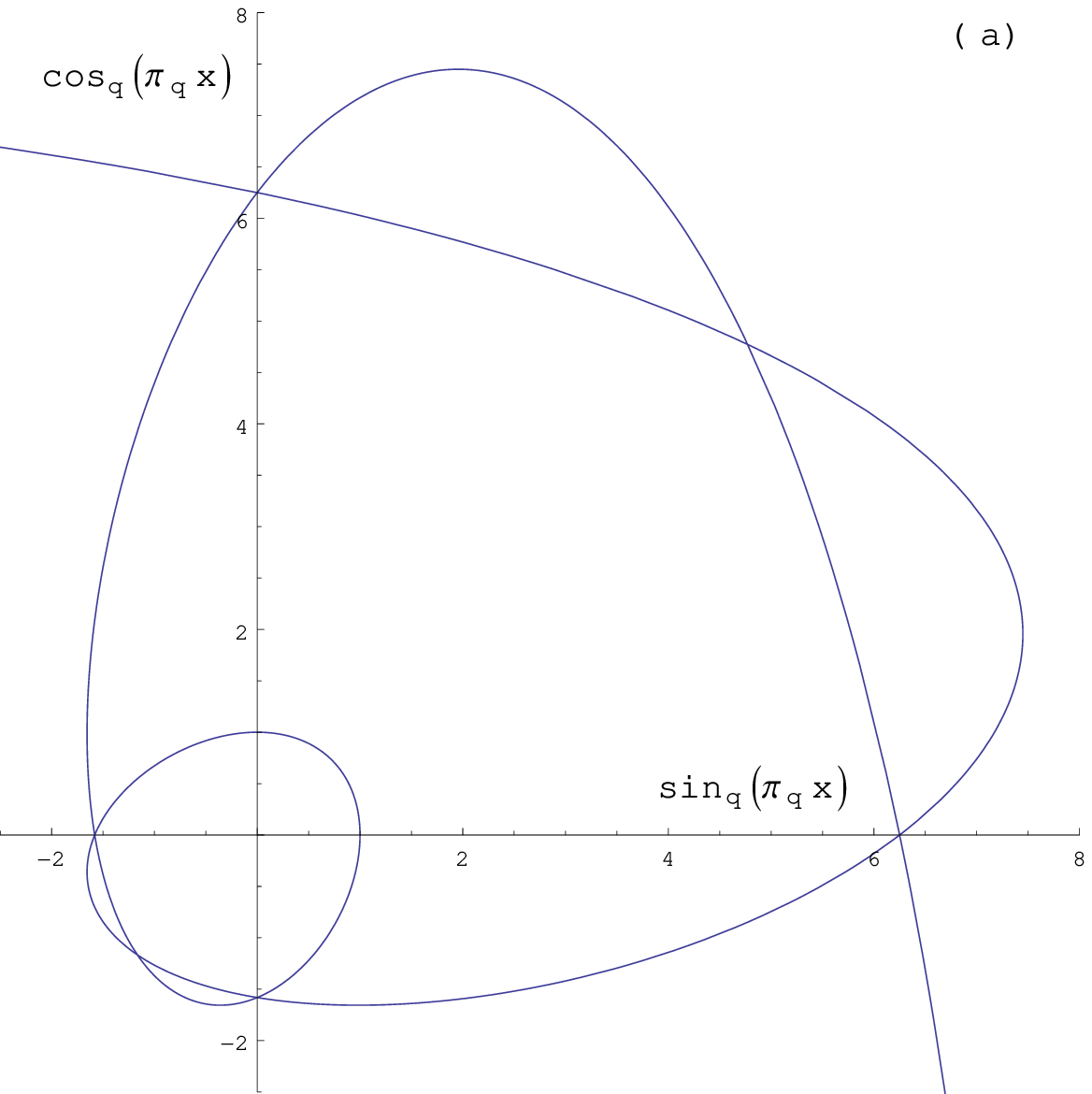}
\includegraphics[scale=0.60]{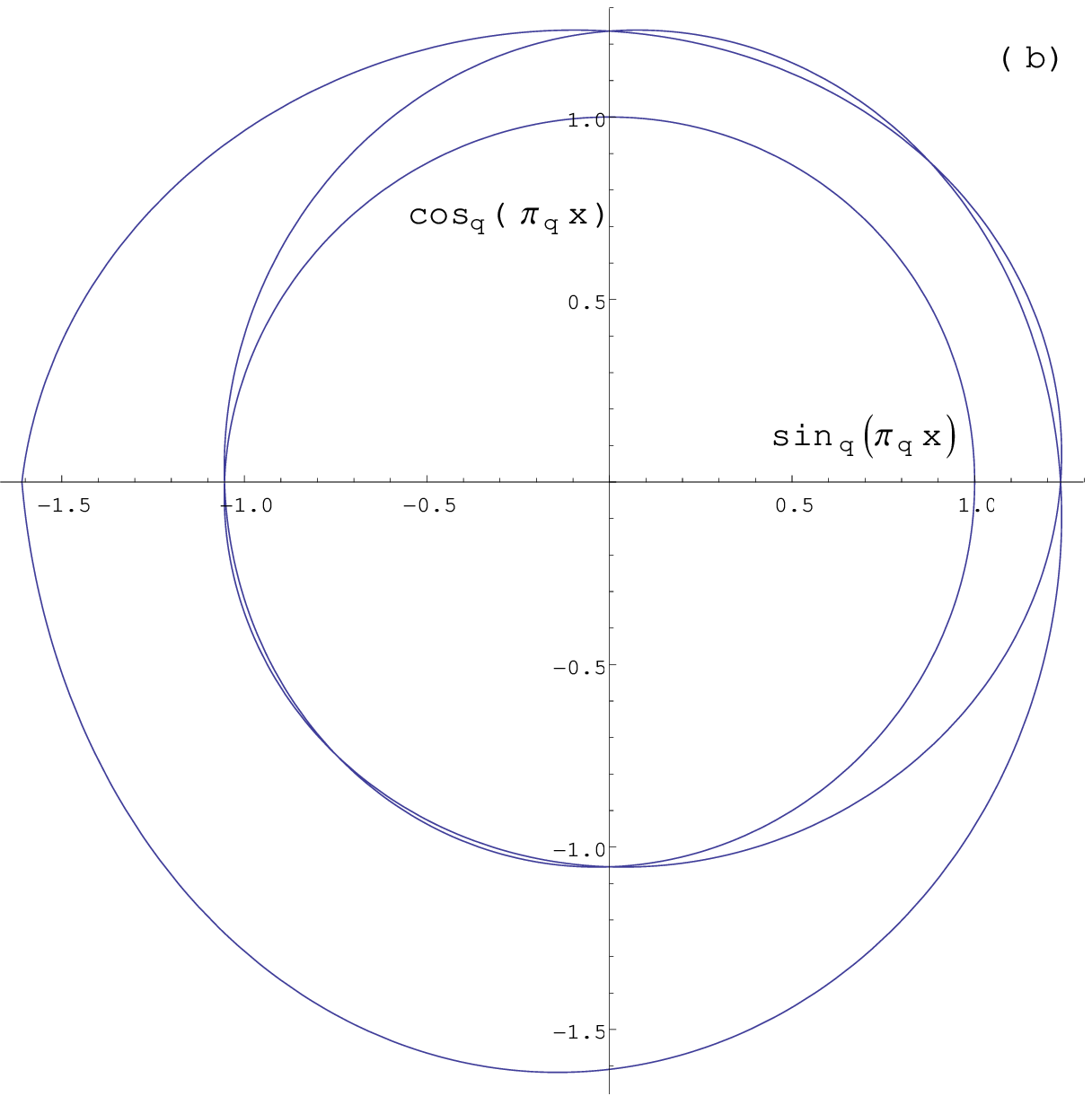}
\end{center}
\caption{\label{fig2} The parametric plot of cosinus-$q$ and sinus-$q$ functions for $q=0.4$ (Fig. \ref{fig2}(a)) and $q=0.9$ (Fig. \ref{fig2}(b)). Notice that for $q$ closer to one the parametric plots is closer to circle as is for the (standard) cosinus and sinus functions.}
\end{figure}
We can summarize the results obtained so far by noting that the thread we followed has been that of exploiting the umbral correspondence between $q$-Gaussian and $q$-Lorentzian functions to infer the existence of $\pi_{q}$, which has been then used to infer the $q$-Wallis formula and a family of trigonometric like functions defined in terms of $q$-Euler infinite product. 

\section{Umbral correspondence and $q$-Bessel function {and $q$-Hermite polynomials}}

The functions ${_{q}e}(x)$ and $C_{0}^{(q, 1)}(x)$ have been shown to be useful study cases. We remind that the first is an eigenfunction of the Jackson derivative \cite{FHJacson1909, GDattoli96, GDattoli97}
\begin{equation}\label{eq29}
{_{q}D_{x}} = \frac{1-q^{x \partial_{x}}}{(1-q)x}.
\end{equation}
By recalling indeed that 
$q^{x \partial_{x}} f(x) = f(q x)$ \cite{GDattoli17a, dzielo}
we find
\begin{equation}\label{eq31}
{_{q}D_{x}} f(x) = \frac{f(x) - f(qx)}{(1-q) x}.
\end{equation}
Thus, getting
\begin{equation}\label{eq32}
{_{q}D_{x}} x^{n} = [n]_{q}\, x^{n-1},
\end{equation}
which yields
\begin{align}\label{eq33}
\begin{split}
{_{q}D_{x}}\, {_{q}e}(\lambda x) & = {_{q}D_{x}}\, \sum_{r=0}^{\infty} \frac{(-\lambda x)^{r}}{[r]_{q}!} = \sum_{r=0}^{\infty} (-\lambda)^{r} {_{q}D_{x}} \frac{x^{r}}{[r]_{q}!} = \sum_{r=1}^{\infty} [r]_{q} \frac{(-\lambda x)^{r-1}}{[r]_{q}!} \\
& = -\lambda\, {_{q}e(\lambda x)}
\end{split}
\end{align}
It is therefore easy to infer that the $(q,1)$-Tricomi satisfies the eigenvalue equation
\begin{equation}\label{eq34}
\partial_{x} x {_{q}D_{x}} C_{0}^{(q, 1)}(\lambda x) = - \lambda C_{0}^{(q, 1)}(\lambda x).
\end{equation}
Regarding the function
\begin{equation}\label{eq35}
C_{0}^{(q, q)}(x) = \sum_{r=0}^{\infty}\frac{(-x)^{r}}{([r]_{q}!)^{2}}
\end{equation}
we argue that
\begin{equation}\label{eq36}
{_{q}D_{x}} x\, {_{q}D_{x}} C_{0}^{(q, q)}(\lambda x) = - \lambda C_{0}^{(q, q)}(x)
\end{equation}
The procedure we have envisaged opens the pathway not only for the study of $q$- Bessel like functions but also for family of special $q$-polynomials.

The properties of the Jackson derivative can e.g. be exploited to derive families of two variables Hermite polynomials, we find e.g.
\begin{equation}\label{eq37}
\E^{y\, {_{q}D_{x}^{2}}} x^{n} = \sum_{r=0}^{\lfloor n/2 \rfloor} \frac{[n]_{q}! y^{r} x^{n-2r}}{[n-2r]_{q}! \, r!} = H_{n}^{(q, 1)}(x, y),
\end{equation}
which are shown to satisfy the recurrences
\begin{align}\label{eq38}
\begin{split}
{_{q}D_{x}} H_{n}^{(q, 1)}(x, y) & = [n]_{q} H_{n-1}^{(q, 1)}(x, y) \\
\partial_{y} H_{n}^{(q, 1)}(x, y) & = [n]_{q} [n-1]_{q} H_{n}^{(q, 1)}(x, y)
\end{split}
\end{align}
and are the natural solution of the $q$-heat equation
\begin{equation}\label{eq39}
\partial_{y} F(x, y) = {_{q}D_{x}^{2}} F(x, y), \qquad F(x, 0) = x^{n}.
\end{equation}
It is evident that $(q,1)$-polynomials have a twofold nature, as summarized by the associated generating function
\begin{equation}\label{eq40}
\sum_{n=0}^{\infty} \frac{t^{n}}{[n_{q}]!} H_{n}^{(q, 1)}(x, y) = \E^{y t^{2}} {_{q}e}(xt).
\end{equation}
We are therefore allowed to keep the ordinary derivative with respect to $x$, which yields
\begin{align}\label{eq41}
\begin{split}
\partial_{x} H_{n}^{(q, 1)}(x, y) &= \sum_{r=0}^{\lfloor n/2 \rfloor} \frac{[n]_{q}!\, y^{r} (n-2r) x^{n-2r-1}}{[n-2r]_{q}! r!} \\
& = \frac{n}{x} H_{n}^{(q, 1)}(x, y) - \frac{2y}{x} [n]_{q} [n-1]_{q} H_{n-2}^{(q, 1)}(x, y).
\end{split}
\end{align}
The previous recurrence, along with those reported in \eqref{eq38}, can eventually be exploited to end up with the mixed differential equation
\begin{equation}\label{eq42}
x \partial_{x}\, H_{n}^{(q, 1)}(x, y) + 2y {_{q}D_{x}}^{2} H_{n}^{(q, 1)}(x, y) = n H_{n}^{(q, 1)}(x, y).
\end{equation}
The same procedure can be exploited to infer the existence of other families of two variables Hermite $(1, q)$ and $(q, q)$ which will be discussed in a forthcoming investigation.

\section{Conclusion}

In this paper we have shown that the $q$-calculus for $0 < q < 1$ can be merged with the umbral calculus. We have also argued that the properties of $q$-functions can be treated quite straightforwardly by using the corresponding umbral images.

In addition, it should be also mentioned that in the physics literatures \cite{ABudini15, CTsallis09} exist another type of the the $q$-exponential function which was introduced by C. Tsallits
\begin{equation}\label{27/03-1}
\tilde{e}_{q}(x) = [1 + (1-q) x]^{1/(1-q)}, \qquad \lim_{q\to 1^{-}} \tilde{e}_{q}(x) = \exp(x).
\end{equation}
and used for the definition of $q$-Gaussian $\tilde{e}_{q}(-x^{2})$. They have been exploited for describing, e.g.,  the normal and tumoral melanocytes \cite{PCAdaSivla}, the momentum distribution of cold atoms in dissipative optical lattices \cite{PDouglas06}, and so on. The $q$-Gaussian umbral form reads
\begin{equation}\label{28/03-1}
\tilde{e}_{1-Q}(-x^{2}) = \E^{- Q x^{2} \hat{d}_{z}} \psi_{0}, \qquad \text{where} \quad \hat{d}_{z}^{\mu}\psi_{0} = \frac{\Gamma(1+1/Q)}{\Gamma(1 + 1/Q - \mu)},
\end{equation} 
where $q = 1-Q$, and it could open the new possibilitless of calculation 
obtained by employing numerical computations. Using the umbral representation of the Tsallis Gaussian and the properties of the ordinary Gaussian integrals 
we obtain for example
\begin{align}\label{4/04-1}
\begin{split}
\int_{-\infty}^{\infty} \tilde{e}_{1-Q}(x^{2}) \D x & = \int_{-\infty}^{\infty} \E^{- Q x^{2} \hat{d}_{z}} \psi_{0} \D x = \int_{-\infty}^{\infty} \E^{-Q x^{2} \hat{d}_{z}} \D x\, \psi_{0} \\
& = \sqrt{\pi/Q}\, \hat{d}_{z}^{-1/2} \psi_{0} = \sqrt{\pi/Q}\, \frac{\Gamma(1+1/Q)}{\Gamma(3/2 + 1/Q)}
\end{split}
\end{align}
and from Eq. \eqref{3/04-1} we introduce the new kind of Hermite polynomials with generating function
\begin{equation}\label{4/04-2}
 \E^{y x} \tilde{e}_{1-Q}(x^{2}) = \E^{y x - Q x^{2} \hat{d}_{z}} \psi_{0} = \sum_{n=0}^{\infty} \frac{x^{n}}{n!} H_{n}(y, -Q \hat{d}_{z}) \psi_{0}.
\end{equation}
This last point opens the possibility of developing new speculations on the polynomials associated with Tsallis type distributions.

In this paper we have shown that umbral methods allow to treat different forms of $q$-calculus from the same point of view. The Tsallis exponential is indeed different from the cases treated in the previous sections. In a forthcoming investigation we will see how the procedure we have defined in these concluding comments is tailor suited to get further progress in the mathematical handling of the Tsallis distributions.

\section*{Acknowledgents}

K.G. was supported by the Polish National Center for Science (NCN) research grant OPUS12 no. UMO-2016/23/B/ST3/01714 and the Polish National Agency for Academic Exchange (NAWA) Programme im. Bekker, project no. PPN/BEK/2018/1/00184.

K.G. would like to thank the ENEA Frascati for warm hospitality.

\end{document}